\documentclass[10pt]{amsart}
\usepackage{amsmath,amsfonts,amssymb,color,a4,epsfig,graphics}
\usepackage[english]{babel}
\usepackage{enumerate}
\usepackage{xcolor}
\RequirePackage{ifthen}
\usepackage{tikz}
\def\build#1_#2^#3{\mathrel{\mathop{\kern 0pt#1}\limits_{#2}^{#3}}}
\usetikzlibrary{arrows}

\newcommand{\ran}{{\rm ran }}
\newcommand{\R}{{\mathbb{R}}}
\newcommand{\C}{{\mathbb{C}}}
\newcommand{\Z}{{\mathbb{Z}}}
\newcommand{\N}{{\mathbb{N}}}

\newcommand{\Cc}{\mathcal{C}}
\newcommand{\Dc}{\mathcal{D}}
\newcommand{\Ec}{\mathcal{E}}
\newcommand{\Gc}{\mathcal{G}}
\newcommand{\Ic}{\mathcal{I}}
\newcommand{\Nc}{\mathcal{N}}

\newcommand{\Vc}{\mathcal{V}}
\def\Hol{{\rm Hol}}

\def\rmi{{\rm i}}

\newtheorem{theorem}{Theorem}[section]
\newtheorem{proposition}[theorem]{Proposition}
\newtheorem{lemma}[theorem]{Lemma}
\newtheorem{remark}[theorem]{Remark}
\newtheorem{example}[theorem]{Example}
\newtheorem{definition}[theorem]{Definition}
\newtheorem{corollary}[theorem]{Corollary}

\begin{document}

\title[The  magnetic Laplacian acting
  on discrete cusps]{The  magnetic Laplacian acting
  on discrete cusps}  
\author{Sylvain Gol\'enia}
\address{Institut de Math\'ematiques de
    Bordeaux,
351, cours de la Lib\'eration
F-33405 Talence cedex}    
\email{\tt sylvain.golenia@math.u-bordeaux.fr}
\author{Fran{\c c}oise Truc}
\address{Grenoble University, Institut Fourier,
Unit{\'e} mixte
 de recherche CNRS-UJF 5582,
 BP 74, 38402-Saint Martin d'H\`eres Cedex, France}
\email{\tt francoise.truc@ujf-grenoble.fr}
\subjclass[2010]{34L20, 47A10, 05C63, 47B25, 47A63, 81Q10}
\keywords{discrete magnetic Laplacian, locally finite graphs,
  eigenvalues, asymptotic, form-domain} 
\begin{abstract}
We introduce the notion of discrete cusp for a weighted graph. In this context, we prove
that the form-domain of the magnetic Laplacian and that of the
non-magnetic Laplacian can be different. We  
establish the emptiness of the essential spectrum and compute the
asymptotic of eigenvalues for the magnetic Laplacian.  
\end{abstract}

\maketitle

\section{Introduction}

The spectral theory of discrete Laplacians on graphs has drawn a lot
of attention for decades. The spectral analysis of the Laplacian
associated to a graph is strongly related to the geometry of the
graph. Moreover, graphs are discretized versions of manifolds. In 
\cite{AbTr, GoMo}, it is shown that for a manifold with cusps, adding
a magnetic field can drastically destroy the essential spectrum of
the Laplacian. The aim of this article is to go along this line in a
discrete setting.

We recall some standard definitions of graph theory. A \emph{graph} is
a triple
$\Gc:=(\mathcal{E},\mathcal{V}, m)$, where $\mathcal{V}$ is a countable
set (the \emph{vertices}), $\mathcal{E}:\mathcal{V} \times
\mathcal{V}\rightarrow
\mathbb{R}_{+}$ 
is symmetric, and $m:\Vc\to (0,
\infty)$ is a weight. 
We
say that $\Gc$ is \emph{simple} if $m=1$ and $\Ec:\Vc\times \Vc\to
\{0,1\}$. 

Given  $x,y\in \mathcal{V}$, we say that $(x,y)$ is an edge (or $x$ and $y$ are
\emph{neighbors}) if $\mathcal{E}(x,y)>0$. We denote this
relationship by $x\sim y$ and the set of neighbors of $x$ by
$N_\Gc(x)$.
We say that there is a \emph{loop} at $x\in\mathcal{V}$
if $\mathcal{E}(x,x)>0$.   A graph is \emph{connected} if for all
$x,y\in\mathcal{V}$, there exists a path $\gamma$ joining $x$ and
$y$. Here, $\gamma$ is a sequence
$x_{0},~x_{1},...,x_{n}\in\mathcal{V}$ such that $x=x_{0},~y=x_{n}$,
and $x_{j}\sim x_{j+1}$ for all $0\leq j \leq n-1$. In this case, we
set $|\gamma|:=n$. 
  A graph $\Gc$ is \emph{locally finite} if $|N_\Gc(x)|$ is finite for
all $x\in\mathcal{V}$. In the sequel, we assume that:
\begin{center}
{\bf All graphs are locally finite, connected with no loops.}
\end{center}
We endow a graph $\Gc:=(\Ec, \Vc, m)$ with the metric
$\rho_{\mathcal{G}}$ defined by 
\begin{align*}
\rho_{\Gc}(x,y):=\inf
\{|\gamma|,~\gamma~\mbox{is~a~path~joining~$x$~and~$y$}\}.
\end{align*}
The space of complex-valued functions acting on the set of vertices $\Vc$ is
denoted by $C(\Vc):=\lbrace f:\Vc \to {\C}\rbrace$. Moreover,
 $C_c(\Vc)$ is the subspace of $C(\Vc)$  of functions with finite
 support. We consider the Hilbert space
\[\ell^2(\Vc, m):=\left\{ f\in C(\Vc), \sum_{x\in
\Vc} m(x)|f(x)|^2 <\infty \right\}\]
with the scalar product $\langle f,g\rangle:= \sum_{x\in \Vc}m(x)
\overline{f(x)}g(x)$. 

We equip $\Gc$  with a \emph{magnetic potential} $\theta:\Vc\times
\Vc\to  {\R/2\pi \Z}$ such that we have
$\theta_{x,y}:=\theta(x,y)= -\theta_{y,x}$ and $\theta(x,y):=0$ if
$\Ec(x,y)=0$. We define the Hermitian form 
\[ Q_{\Gc,\theta} (f):= \frac{1}{2}\sum_{x,y \in \Vc} \Ec(x,y) \left| f(x) -e^{\rmi
  \theta_{x,y}}f(y)\right|^2,\]
 for all $f\in \Cc_c(\Vc)$. 
The associated \emph{magnetic
  Laplacian}  is the unique non-negative self-adjoint
operator $\Delta_{\Gc,
  \theta}$  satisfying $\langle f, \Delta_{\Gc,
  \theta} f\rangle_{\ell^2(\Vc,m)}= Q_{\Gc, \theta} (f)$, for all $f\in
\Cc_c(\Vc)$.  It is the \emph{Friedrichs extension} of
$\Delta_{\Gc, \theta}|_{\Cc_c(\Vc)}$, e.g., \cite{CTT3, RS},
 where 
\[(\Delta_{\Gc, \theta} f)(x)=\frac{1}{m(x)}\sum_{y\in \Vc}\Ec(x,y) \left(f(x)- e^{\rmi \theta_{x,y}}f(y)\right),\]
for all $f\in \Cc_c(\Vc)$. 
 We set 
\[\deg_\Gc(x):=\frac{1}{m(x)}\sum_{y\in \Vc} \Ec(x,y),\]
the \emph{degree} of $x\in \Vc$. 
We see easily that $\Delta_{\Gc, \theta}\leq 2 \deg_\Gc(\cdot)$
in the \emph{form sense}, i.e.,
\begin{align}\label{e:inedeg}
0\leq\left\langle f,\Delta_{\Gc, \theta} f\right\rangle
\leq  \langle f, 2\deg_{\Gc}(\cdot) f\rangle, \mbox{ for all }f\in \Cc_c(\Vc).
\end{align}
Moreover,  setting $\tilde \delta_x(y):=  m^{-1/2}(x)\delta_{x,y}$ for any $x,y \in \Vc$, $\langle \tilde \delta_x, \Delta_{\Gc, \theta} \tilde \delta_x\rangle = \deg_\Gc(x)$,
so $\Delta_{\Gc, \theta}$ is bounded if and only if $\sup_{x\in
  \Vc}\deg_\Gc(x)$ is finite, e.g.\ \cite{KL, Go}.

Another consequence of \eqref{e:inedeg} is 
\begin{align}\label{e:inclusion}
\Dc\left(\deg_{\Gc}^{1/2}(\cdot)\right)\subset \Dc\left(\Delta_{\Gc,
    \theta}^{1/2}\right),
\end{align}
where $\Dc\left(\deg_{\Gc}^{1/2}(\cdot)\right):= \left\{f\in \ell^2(\Vc, m), \deg_\Gc(\cdot)f \in  \ell^2(\Vc, m)\right\}$.
However, the equality of the form-domains
\begin{align}\label{e:equa}
\Dc\left(\deg_{\Gc}^{1/2}(\cdot)\right)= \Dc\left(\Delta_{\Gc, \theta}^{1/2}\right)
\end{align}
is wrong in general for a simple graph, see \cite{Go, BGK}. In fact if
$\theta=0$, \eqref{e:inclusion} is equivalent to a sparseness
condition and holds true for planar simple 
graphs, see \cite{BGK}. We refer to \cite{BGKLM} for a magnetic
sparseness condition.  On a general weighted graph, if \eqref{e:equa}
holds true, 
\[\sigma_{\rm ess}(\Delta_{\Gc, \theta})=\emptyset \Leftrightarrow (\Delta_{\Gc, \theta}+1)^{-1} \mbox{ is compact} \Leftrightarrow \lim_{|x|\to \infty} \deg_\Gc(x)=\infty,\]
where $|x|:=\rho_{\Gc}(x_0,x)$ for a given $x_0\in \Vc$. Note that the limit is independent of the choice of $x_0$. 
Besides if the latter is true and if the graph is sparse (simple and planar for instance), \cite{BGK} ensures the following asymptotic of eigenvalues,
\begin{align}\label{e:asympgene}
\lim_{n\to \infty}\frac{\lambda_n \left(\Delta_{\Gc,
  \theta}\right)}{\lambda_n \left(\deg_\Gc(\cdot)\right)}=1,
\end{align}
 where $\lambda_n(H)$ denotes the $n$-th
eigenvalue,  counted with multiplicity, of a self-adjoint operator
$H$, which is bounded from below. 

The technique used in \cite{BGK} does not apply when the graph is a discrete
cusp (thin at infinity),  see Definition \ref{d:dc}. The aim of this article is to establish
new behaviors for the asymptotic of eigenvalues for the magnetic Laplacian in that case, 
and also to prove that the form-domain of the non-magnetic
Laplacian can be different from that of the magnetic Laplacian, see
Theorem \ref{t:main_cp}. We found 
the inspiration by mimicking the continuous case,  which was studied
in \cite{AbTr,GoMo}. 

Let us present a flavour of our results (in particular of Theorem \ref{t:main_cp}) by introducing the following specific example of \emph{discrete cusp} :

\begin{example}\label{t:intro}
Let $n\geq 3$ be an integer and consider $\Gc_1:=(\Ec_1, \Vc_1, m_1)$, where 
\[\Vc_1:=\N, \quad m_1(n):=\exp(-n), \mbox{ and } \Ec_1(n, n+1):= \exp(-(2n+1)/2),\]
for all $n\in \N$ and $\Gc_2:=(\Ec_2, \Vc_2, 1)$ a simple connected finite graph such that $|\Vc_2|=n$.
Set $\theta_1:=0$ and $\theta_2$ such that $\Hol_{\theta_2}\neq 0$.
 Let  $\Gc:=(\Ec, \Vc, m)$ be the \emph{twisted  Cartesian product} $\Gc_1\times_{\Vc_2} \Gc_2$, given by:
\begin{align*}
\left\{\begin{array}{rl}
m(x,y):=& m_1(x),
\\
\Ec\left((x,y),(x',y')\right):= 
&\Ec_1(x,x')\times \delta_{y, y'} +
\delta_{x,x' } \times \Ec_2(y,y'),
\\
\theta\left((x,y),(x',y')\right):= 
&
\delta_{x,x' } \times \theta_2(y,y'),
\end{array}\right.
\end{align*}
for all $x,x'\in \Vc_1$ and $y,y'\in \Vc_2$. 
Then there exists
  a constant $\nu>0$ such that for all $\kappa\in \R/\nu\Z$
\begin{align*}
\sigma_{\rm ess}(\Delta_{\Gc, \kappa \theta})=\emptyset \Leftrightarrow
\Dc\left(\Delta_{\Gc, \kappa \theta}^{1/2}\right) =
\Dc\left(\deg_{\Gc}^{1/2}(\cdot)\right)\Leftrightarrow \kappa \neq 0
\mbox{  in } \R/\nu\Z 
\end{align*}
Moreover: 
\begin{enumerate}[1)]
\item When  $\kappa\neq 0 \mbox{  in } \R/\nu \Z$, we have:
\[\lim_{\lambda \to \infty}\frac{\Nc_\lambda\left(\Delta_{\Gc, \kappa
    \theta}\right)}{\Nc_\lambda\left(\deg_\Gc(\cdot)\right)}=1,\]
where $\Nc_\lambda(H):= \dim \ran 1_{]-\infty, \lambda]}(H)$ for a
self-adjoint operator $H$. 
\item When $\kappa =0$ in $\R/\nu \Z$, the absolutely continuous part
of the $\Delta_{\Gc, \kappa     \theta}$ is
\[\sigma_{\rm ac} \left(\Delta_{\Gc, \kappa     \theta}\right)= \left[e^{1/2}+e^{-1/2}-2, e^{1/2}+e^{-1/2} 
+2\right],\]
with multiplicity $1$ and 
\[\lim_{\lambda \to \infty}\frac{\Nc_\lambda\left(\Delta_{\Gc, \kappa
    \theta}P^\perp_{{\rm
    ac}, \kappa}\right)}{\Nc_\lambda\left(\deg_\Gc(\cdot)\right)}=\frac{n -1}{n},\]
where $P_{{\rm ac}, \kappa }$ denotes the projection onto the a.c.\ part of
$\Delta_{\Gc, \kappa     \theta}$. 
\end{enumerate}
\end{example}

We now describe heuristically the phenomenon. Compared with the first
case, the constant $(n-1)/n$ that appears in the second case encodes
the fact that a part of the wave packet diffuses. Moreover, switching on the magnetic field is not a gentle perturbation 
because the form domain of the operator is changed.

By Riemann-Lebesgue Theorem, the particle, which is localized in the
a.c.\ part of the 
operator, escapes from every compact set. More precisely, for
a finite subset $X\subset \Vc$ and all $f\in \Dc(\Delta_{\Gc, 0})$ 
\[\|1_{\rm X}(\cdot)  e^{\rmi t \Delta_{\Gc, 0}}P_{{\rm ac}, 0} f\| \to 0, \mbox{ as }
t\to \infty.\]
In the first case, when the magnetic potential is active, the spectrum of
$\Delta_{\Gc, \kappa\theta}$ is purely discrete. The particle cannot
diffuse anymore.  More precisely, for
a finite subset $X\subset \Vc$ and an eigenvalue $f$
of $\Delta_{\Gc, \kappa\theta}$ such that $f|_X \neq 0$, there is $c>0$ such that:
\[\frac{1}{T}\int_0^T \|1_{\rm X}(\cdot)e^{\rmi t \Delta_{\Gc, \kappa
    \theta}}f\|^2\, dt \to c, \mbox{ as } T \to \infty.\]
The particle is \emph{trapped} by the magnetic field.

\begin{align*}
\begin{tikzpicture}[scale=0.55]
\def\test{0.75}
\def\xmax{10}
  \foreach \x in {1,2,...,\xmax}
{\fill[color=black]({(\x)/\test}, 0)circle(.3mm);
\fill[color=black]({(\x)/\test+.25/\test/\x}, 4*0.75/\x)circle(.3mm);
\fill[color=black]({(\x)/\test-1/\test/\x}, 4*1/\x)circle(.3mm);
\draw(\x/\test, 0)--(\x/\test+.25/\test/\x, 4*0.75/\x);
\draw(\x/\test, 0)--(\x/\test-1/\test/\x, 4*1/\x);
\draw(\x/\test+.25/\test/\x, 4*0.75/\x)--(\x/\test-1/\test/\x, 4*1/\x);
\draw(\x/\test, 0)--({(\x+1)/\test}, 0);
\draw(\x/\test-1/\test/\x, 4*1/\x)--({(\x+1)/\test-1/\test/(\x+1)}, {4*1/(\x+1)});
\draw(\x/\test+.25/\test/\x, 4*0.75/\x)--({(\x+1)/\test+.25/\test/(\x+1)}, {4*.75/(\x+1)});
}
\fill[color=black]({(\xmax+1)/\test}, 0)circle(0.3mm);
\fill[color=black]({(\xmax+1)/\test+.25/\test/(\xmax+1)}, {4*0.75/(\xmax+1)})circle(.3mm);
\fill[color=black]({(\xmax+1)/\test-1/\test/(\xmax+1)}, {4*1/(\xmax+1)})circle(.3mm);
\draw({(\xmax+1)/\test}, 0)--({(\xmax+1)/\test+.25/\test/(\xmax+1)}, {4*0.75/(\xmax+1)});
\draw({(\xmax+1)/\test}, 0)--({(\xmax+1)/\test-1/\test/(\xmax+1)}, {4*1/(\xmax+1)});
\draw({(\xmax+1)/\test+.25/\test/(\xmax+1)}, {4*0.75/(\xmax+1)})--({(\xmax+1)/\test-1/\test/(\xmax+1)}, {4*1/(\xmax+1)});
%
%
\path({(\xmax+1)/\test+1}, 0) node {$\cdots$};
\draw[red, thick, ->]   (9,-0.5) .. controls (10,-0.5) and (10,1-0.5)
.. (9,2-0.5);
.. (-7,2-0.5);
\draw[blue, thick, ->]   (2,4) .. controls (4,2.5) 
.. (6,2);
.. (-6,2);
\path[blue](4,4) node {Diffusion};
\path[red](11,2) node {Magnetic effect};
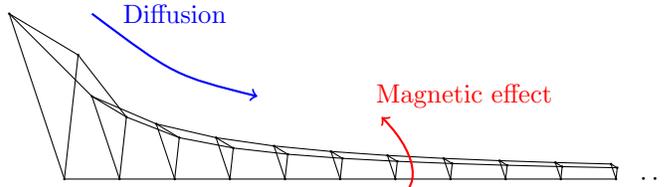
\path(8, -1) node {\emph{Representation of  a discrete cusp:}};
\path(8, -2) node {\emph{The magnetic field traps the particle by
  spinning it,}};
\path(8, -3) node {\emph{whereas its absence lets the particle diffuse.}};
\end{tikzpicture}
\end{align*}

We now describe the structure of the paper. In Section \ref{s:holo}, we
recall some properties of the holonomy of a magnetic
potential. In Section \ref{s:setting} we present 
 our main hypotheses and several notions of (weighted)
product for graphs. We introduce the notion of discrete cusp and
analyze it under the light of the radius of injectivity.
 Then
in Section \ref{s:destroy} we give a criteria concerning the absence of essential
spectrum. Next, in Section \ref{s:asymp}, we refine the 
analysis and give our central theorem, a general statement for discrete cusps, computing the form domain and  the asymptotic of eigenvalues. We 
finish the  
section by proving Theorem \ref{t:intro}.

\noindent\textbf{Notation:} $\N$ denotes the set of non negative
integers and $\N^*$ that of the positive integers. We denote by
$\Dc(H)$ the domain of an operator $H$. Its (essential) spectrum is
denoted by $\sigma(H)$ (by $\sigma_{\rm ess}(H)$). We set $\delta_{x,y}$ equals $1$ if and only if $x=y$ and $0$ otherwise and given a set $X$, $1_{X}(x)$ equals $1$ if $x\in X$ and $0$ otherwise.

\noindent\textbf{Acknowledgments:}
We would like to thank Colette Ann\'e, Michel Bonnefont, Yves Colin de Verdi\`ere, Matthias Keller, and Sergiu Moroianu for useful discussions. SG and FT were partially supported by the ANR project GeRaSic (ANR-13-BS01-0007-01) and by SQFT (ANR-12-JS01-0008-01).

\section{Main results}
\subsection{Holonomy of a magnetic potential}\label{s:holo}
We recall some facts about the gauge theory of magnetic
fields, see \cite{CTT3, HiSh} for more details and also \cite{LLPP} for a different point of view. We recall that a \emph{gauge
  transform} $U$ is the unitary map on 
$\ell^2(\Vc, m)$
defined by
\[(Uf)(x)= u_x f(x),\]
where $(u_x)_{x\in \Vc}$
is  a sequence of complex numbers
  with $|u_x|\equiv 1 $ (we write $u_x=e^{\rmi \sigma_x}$).
The map $U $ acts on the quadratic forms  $Q_{\Gc,\theta}$
by $U^\star (Q_{\Gc, \theta})(f) = Q_{\Gc, \theta} (Uf)$, for all $f\in \Cc_c(\Vc)$.
The magnetic potential  $U^\star (\theta) $ is defined by:
\[U^\star (Q_{\Gc, \theta})=Q_{\Gc, U^\star (\theta)}.\] 
More explicitly, we get:
\[U^\star (\theta)_{xy}=  \theta_{x,y}+\sigma_y
  -\sigma_x. \]
We turn to the definition of the flux of a magnetic potential, the
\emph{Holonomy}. 

\begin{proposition}\label{p:Z1}
Let us denote by $Z_1(\Gc)$ the
space of cycles of $\Gc$. It is is a
free $\Z-$module with a basis of
geometric cycles 
$\gamma =(x_0, x_1)+(x_1,x_2)+ \ldots  + (x_{N-1} ,x_N )$
 with, for $i=0, \cdots, N-1$,
$\Ec(x_i ,x_{i+1})\neq 0$, and $x_N=x_0$.
We define the \emph{holonomy map} 
$\Hol_\theta :Z_1(\Gc)\to \R/2\pi \Z $, 
by 
\[\Hol_\theta \left((x_0, x_1)+(x_1,x_2)+ \cdots  + (x_N ,x_0 ) \right):=
\theta _{x_0,x_1}+ \cdots +\theta_{x_N, x_0}.\]  
Then
\begin{enumerate}[1)]
  \item The map $\theta \mapsto \Hol_\theta $ is surjective
onto  ${\rm Hom }_\Z(Z_1(\Gc),\R/2 \pi \Z )$.
  \item $\Hol_{\theta_1} =\Hol_{\theta_2}$ if and only if there exists
    a gauge transform $U$ 
so that $U^\star (\theta_2 )= \theta_1$.
\end{enumerate}
In consequence  $\Hol_{\theta_1} =\Hol_{\theta_2}$ if and only if the magnetic Laplacians $\Delta_{\Gc, \theta_1} $ and $\Delta_{\Gc, \theta_2} $ 
are unitarily equivalent.
\end{proposition}
\begin{lemma}\label{l:keylemm}
Let $\Gc:=(\Ec, \Vc, m)$ be a connected graph such that $1\in\ker  
\Delta_{\Gc,0}$. Let $\theta$ be magnetic
potential. Then  $\ker  
\Delta_{\Gc,\theta} \neq \{0\}$ if and only if $\Hol_\theta =0$. 
\end{lemma}
\begin{remark}
By construction of the Friedrichs extension, the domain of $\Delta_{\Gc,0}$ is given by
\begin{align*}
\Dc(\Delta_{\Gc,0})&=\left\{f\in \ell^2(\Vc, m), x\mapsto \frac{1}{m(x)}\sum_{y\in \Vc} \Ec(x,y) (f(x)-f(y)) \in \ell^2(\Vc, m) \right\} 
\\
& \quad \bigcap \overline{\Cc_c(\Vc)}^{(\|\cdot\|^2 + Q_{\Gc, 0}(\cdot))^{1/2}}.
\end{align*}
The hypothesis $1\in\ker  
\Delta_{\Gc,0}$ is trivially satisfied if $\Gc$ is a finite graph. In general, it is satisfied if and only if:
\begin{itemize}
\item[$(*)$]  $1$ belongs to the closure of $\Cc_c(\Vc)$ with respect to the norm $(\|\cdot\|^2 + Q_{\Gc, 0}(\cdot))^{1/2}$.
\end{itemize}
A sufficient condition to guarantee $(*)$ is that the following two conditions hold true:
\begin{enumerate}[1)]
\item $\Gc$ is of finite volume, i.e.,
such that $\sum_{x\in   \Vc} m(x)<\infty$,
\item $\Delta_{\Gc,0}$ is essentially self-adjoint on $\Cc_c(\Vc)$.
\end{enumerate}
\end{remark}

\proof
If $\Hol_\theta =0$ then $\Delta_{\Gc,\theta}$ is unitarily equivalent
to $\Delta_{\Gc,0}$ by Proposition \ref{p:Z1} and $1\in \ker(\Delta_{\Gc,0})\neq \{0\}$ by hypothesis. 

Conversely, let $f\neq 0 $ with $ \Delta_{\Gc,\theta} f=0$
and hence $Q_{\Gc,\theta}(f)=0.$
This implies that all terms in the expression of
$Q_{\Gc,\theta}(f)$
vanish. In particular, if $\Ec(x,y)\neq 0$ we have 
\begin{equation}\label{e:edge}
f(x) =e^{\rmi
    \theta_{x,y}}f(y).
\end{equation}
Assume that there is a cycle $\gamma =(x_0, x_1,\ldots, x_N=x_0) $,
such that $\Hol_\theta(\gamma)\neq 0$. Using \eqref{e:edge}, 
we obtain that
\[  f(x_i)= e^{-\rmi {\rm  Hol}_\theta (\gamma)} f(x_i)~.\]
for all $i=0, \ldots, N-1$. Therefore $f|_\gamma=0$. 
Then, since $f\neq 0$, there is $x\in \Vc$ such that $f(x)\neq
0$. Using again \eqref{e:edge} and by connectedness between $x$ and
$\gamma$, it yields that $f(x)=0$. Contradiction. Therefore if there
exists $f\in\ker\left(\Delta_{\Gc, \theta}\right) \setminus\{0\}$ then
  $\Hol_\theta=0$. \qed 

We exhibit the following coupling constant effect.

\begin{corollary}\label{c:coupling}
Let $\Gc:=(\Ec, \Vc, m)$ be a connected graph of finite volume, i.e.,
such that $\sum_{x\in   \Vc} m(x)<\infty$ and let $\theta$ be a magnetic
potential such that $\Hol_\theta\neq 0$. Assume that the function $1$ is in  $\ker  
\Delta_{\Gc,\theta} $. Then there is $\nu\in \R$ such that
\[ \ker  
\Delta_{\Gc,\lambda \theta} \neq \{0\} \Leftrightarrow \lambda =0
\mbox{ in } \R/\nu \Z.\]
\end{corollary}
\proof 
Let $\Phi:  (\R, +) \to  ({\rm Hom }_\Z(Z_1(\Gc),\R/2 \pi \Z ), +)$ be defined by
$\Phi (\lambda):= \Hol_{\lambda \theta}$. It is a homomorphism
of group. Hence its kernel is a subgroup of $ (\R, +)$. In particular it  is either dense with respect to the Euclidean norm or equal to $\nu\Z$ for some $\nu\in
\R$, e.g., \cite[Section V.1.1]{Bou}. Suppose by contradiction that the kernel is dense. 
Since for any  cycle $\gamma$ of
$\Gc$,  the map $\lambda \mapsto \Hol_{\lambda
  \theta}(\gamma)$ is continuous from $\R$ to $\R/2\pi \Z$, we infer that $\Hol_{\lambda
  \theta}(\gamma)=0$ for all $\lambda\in \R$. Hence, $\Phi(\lambda)=0$ for all $\lambda \in \R$. This is a contradiction with $\Hol_\theta\neq 0$. We conclude that there is $\nu\in \R$ such that $\ker(\Phi)= \nu\Z$, i.e.,  using Proposition \ref{p:Z1}, that
\[\{\lambda\in \R, \ker  
\Delta_{\Gc,\lambda\theta} \neq \{0\}\}= \{\lambda\in \R,   
\Hol_{\lambda\theta}=0\}= \nu\Z.\]  
 This ends the proof. \qed
\subsection{The setting}\label{s:setting}
Given $\Gc_1:=(\Ec_1, \Vc_1, m_1)$ and
$\Gc_2:=(\Ec_2, \Vc_2, m_2)$, the \emph{Cartesian product of $\Gc_1$ by $\Gc_2$} is defined by 
$\Gc:=(\Ec, \Vc, m)$, where $\Vc:= \Vc_1\times \Vc_2$. 
\begin{align*}
\left\{\begin{array}{rl}
m(x,y):=& m_1(x)\times m_2(y),
\\
\Ec\left((x,y),(x',y')\right):= 
&\Ec_1(x,x')\times \delta_{y, y'}  m_2(y)+
m_1(x)\delta_{x,x' } \times \Ec_2(y,y'),
\\
\theta\left((x,y),(x',y')\right):= 
&\theta_1(x,x')\times \delta_{y, y'}+
\delta_{x,x' } \times \theta_2(y,y'),
\end{array}\right.
\end{align*}
We denote by $\Gc:= \Gc_1\times \Gc_2$. This definition generalizes the unweighted Cartesian product, e.g., \cite{Ha}. It is used in several places in the literature, e.g., \cite{Ch}[Section 2.6] and in \cite{BGKLM} for a generalization. 
\begin{align*}
\begin{tikzpicture}[ scale=1]
\path(-.5, 0) node {$\cdots$};
\path(6.5, 0) node {$\cdots$};
\fill[color=black](0, 0)circle(.5mm);
\fill[color=black](0.5, 0.75)circle(.5mm);
\fill[color=black](-0.5, 1)circle(.5mm);
\fill[color=black](2+0, 0)circle(.5mm);
\fill[color=black](2+0.5, 0.75)circle(.5mm);
\fill[color=black](2+-0.5, 1)circle(.5mm);
\fill[color=black](4+0, 0)circle(.5mm);
\fill[color=black](4+0.5, 0.75)circle(.5mm);
\fill[color=black](4+-0.5, 1)circle(.5mm);
\fill[color=black](6+0, 0)circle(.5mm);
\fill[color=black](6+0.5, 0.75)circle(.5mm);
\fill[color=black](6+-0.5, 1)circle(.5mm);
\draw[-](0, 0)--(2, 0);
\draw[-](2, 0)--(4, 0);
\draw[-](4, 0)--(6, 0);
\draw[-](0+.5, 0.75)--(2+.5, 0.75);
\draw[-](2+.5, 0.75)--(4+.5, 0.75);
\draw[-](4+.5, 0.75)--(6+.5, 0.75);
\draw[-](0-.5, 1)--(2, 1);
\draw[-](2-.5, 1)--(4, 1);
\draw[-](4-.5, 1)--(6-.5, 1);
\draw[-](0, 0)--(.5, 0.75);
\draw[-](0, 0)--(-.5, 1);
\draw[-](.5, 0.75)--(-.5, 1);
\draw[-](2+0, 0)--(2+.5, 0.75);
\draw[-](2+0, 0)--(2+-.5, 1);
\draw[-](2+.5, 0.75)--(2+-.5, 1);
\draw[-](4+0, 0)--(4+.5, 0.75);
\draw[-](4+0, 0)--(4+-.5, 1);
\draw[-](4+.5, 0.75)--(4+-.5, 1);
\draw[-](6+0, 0)--(6+.5, 0.75);
\draw[-](6+0, 0)--(6+-.5, 1);
\draw[-](6+.5, 0.75)--(6+-.5, 1);
\path(+3, -1) node {\emph{The graph of $\Z\times \Z/3\Z$}};
\end{tikzpicture}
\end{align*}
The terminology is motivated by the following decomposition: 
\begin{align*}
\Delta_{\Gc, \theta} = \Delta_{\Gc_1, \theta_1} \otimes 1 +
1 \otimes \Delta_{\Gc_2,\theta_2},
\end{align*}
where $\ell^2(\Vc,m)\simeq\ell^2(\Vc_1, m_1)\otimes \ell^2(\Vc_2, m_2)$. The spectral theory of $\Delta_{\Gc, \theta}$ is well-understood since
\[e^{\rmi t \Delta_{\Gc, \theta}}= e^{\rmi t \Delta_{\Gc_1, \theta_1}}\otimes e^{\rmi t \Delta_{\Gc_2, \theta_2}}, \mbox{ for } t\in \R.\]
 We refer to \cite{RS}[Section VIII.10] for an introduction to
 the tensor product of self-adjoint operators.

In this paper, we are motivated by a geometrical situation. A hyperbolic manifold of finite volume is the union of a compact part 
and of a cusp, e.g., \cite[Theorem 4.5.7]{Th}. The cusp part can be
seen as the product of $(1, \infty)\times M$, where $(M, g_M)$  
is a possibly disconnected Riemannian manifold, endowed with the metric,
\[y^{-1} (dy^2 + g_M).\]
On the cusp part, the infimum of the radius of injectivity is $0$. 

To analyze the Laplacian on this product one separates the variables and obtain a decomposition which is not of the type of a Cartesian product, 
e.g., \cite[Eq.\ (5.22)]{GoMo} for some details. We aim at mimicking this situation and introduce a \emph{modified Cartesian product}. 
Given $\Gc_1:=(\Ec_1, \Vc_1, m_1)$ and
$\Gc_2:=(\Ec_2, \Vc_2, m_2)$ and $\Ic\subset \Vc_2$, we define the \emph{product of $\Gc_1$ by $\Gc_2$ through $\Ic$} by 
$\Gc:=(\Ec, \Vc, m)$, where $\Vc:= \Vc_1\times \Vc_2$ and
\begin{align*}
\left\{\begin{array}{rl}
m(x,y):=& m_1(x)\times m_2(y),
\\
\Ec\left((x,y),(x',y')\right):= 
&\Ec_1(x,x')\times \delta_{y, y'} \left(\sum_{z\in \Ic} \delta_{y,z}\right)+
\delta_{x,x' } \times \Ec_2(y,y'),
\\
\theta\left((x,y),(x',y')\right):= 
&\theta_1(x,x')\times \delta_{y, y'}+
\delta_{x,x' } \times \theta_2(y,y'),
\end{array}\right.
\end{align*}
for all $x, x'\in \Vc_1$ and $y, y'\in \Vc_2$. We denote $\Gc$ by
$\Gc_1\times_\Ic \Gc_2$. If $\Ic$ is empty, the graph is  disconnected
and of no interest for our purpose. If $|\Ic|=1$,
$\Gc_1\times_\Ic \Gc_2$ is the graph $\Gc_1$ decorated by $\Gc_2$, see
\cite{ScAi} for its spectral analysis in the unweighted case. If
$\Ic=\Vc_2$ and $m=1$, we notice that $\Gc_1\times_\Ic \Gc_2 = \Gc_1\times \Gc_2$.

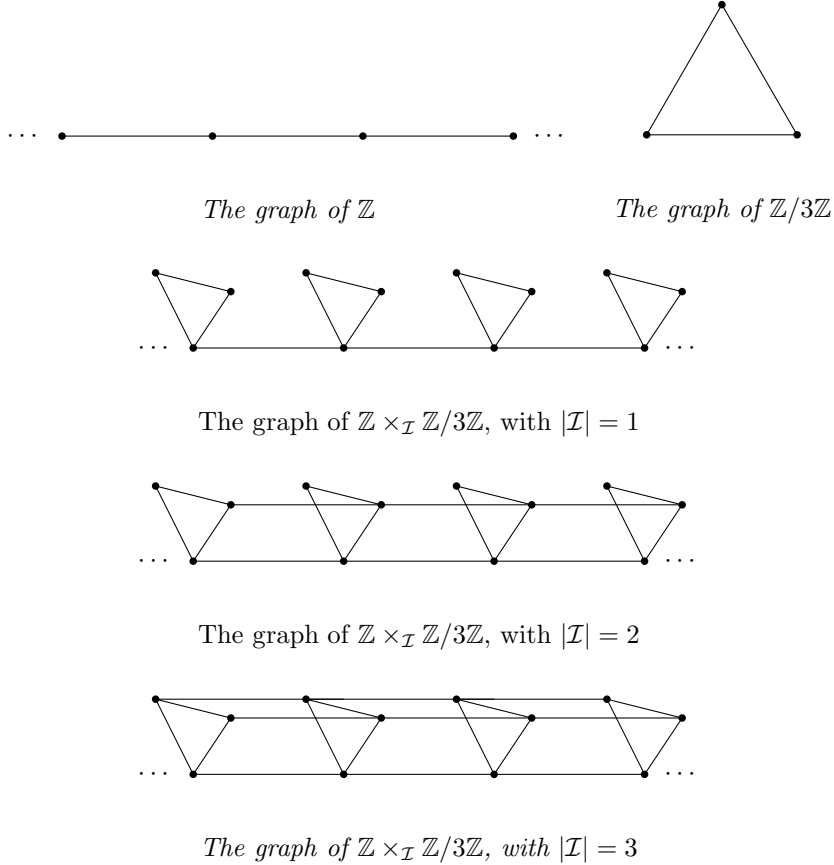
\begin{align*}
\begin{tikzpicture}[ scale=1]
\path(-.5, 0) node {$\cdots$};
\path(6.5, 0) node {$\cdots$};
\fill[color=black](0, 0)circle(.5mm);
\fill[color=black](2+0, 0)circle(.5mm);
\fill[color=black](4+0, 0)circle(.5mm);
\fill[color=black](6+0, 0)circle(.5mm);
\draw[-](0, 0)--(2, 0);
\draw[-](2, 0)--(4, 0);
\draw[-](4, 0)--(6, 0);
\path(+3, -1) node {\emph{The graph of $\Z$}};
\end{tikzpicture}\quad
\begin{tikzpicture}[ scale=1]
\fill[color=black](0, 0)circle(.5mm);
\fill[color=black](1, {sqrt(3)})circle(.5mm);
\fill[color=black](2, 0)circle(.5mm);
\draw[-](0, 0)--(1, {sqrt(3)});
\draw[-](2, 0)--(1, {sqrt(3)});
\draw[-](0, 0)--(2,0);
\path(1, -1) node {\emph{The graph of $\Z/3\Z$}};
\end{tikzpicture}
\end{align*}
\begin{align*}
\begin{tikzpicture}[ scale=1]
\path(-.5, 0) node {$\cdots$};
\path(6.5, 0) node {$\cdots$};
\fill[color=black](0, 0)circle(.5mm);
\fill[color=black](0.5, 0.75)circle(.5mm);
\fill[color=black](-0.5, 1)circle(.5mm);
\fill[color=black](2+0, 0)circle(.5mm);
\fill[color=black](2+0.5, 0.75)circle(.5mm);
\fill[color=black](2+-0.5, 1)circle(.5mm);
\fill[color=black](4+0, 0)circle(.5mm);
\fill[color=black](4+0.5, 0.75)circle(.5mm);
\fill[color=black](4+-0.5, 1)circle(.5mm);
\fill[color=black](6+0, 0)circle(.5mm);
\fill[color=black](6+0.5, 0.75)circle(.5mm);
\fill[color=black](6+-0.5, 1)circle(.5mm);
\draw[-](0, 0)--(2, 0);
\draw[-](2, 0)--(4, 0);
\draw[-](4, 0)--(6, 0);
\draw[-](0, 0)--(.5, 0.75);
\draw[-](0, 0)--(-.5, 1);
\draw[-](.5, 0.75)--(-.5, 1);
\draw[-](2+0, 0)--(2+.5, 0.75);
\draw[-](2+0, 0)--(2+-.5, 1);
\draw[-](2+.5, 0.75)--(2+-.5, 1);
\draw[-](4+0, 0)--(4+.5, 0.75);
\draw[-](4+0, 0)--(4+-.5, 1);
\draw[-](4+.5, 0.75)--(4+-.5, 1);
\draw[-](6+0, 0)--(6+.5, 0.75);
\draw[-](6+0, 0)--(6+-.5, 1);
\draw[-](6+.5, 0.75)--(6+-.5, 1);
\path(+3, -1) node {The graph of $\Z\times_\Ic \Z/3\Z$, with $|\Ic|=1$};
\end{tikzpicture}
\end{align*}
\begin{align*}
\begin{tikzpicture}[ scale=1]
\path(-.5, 0) node {$\cdots$};
\path(6.5, 0) node {$\cdots$};
\fill[color=black](0, 0)circle(.5mm);
\fill[color=black](0.5, 0.75)circle(.5mm);
\fill[color=black](-0.5, 1)circle(.5mm);
\fill[color=black](2+0, 0)circle(.5mm);
\fill[color=black](2+0.5, 0.75)circle(.5mm);
\fill[color=black](2+-0.5, 1)circle(.5mm);
\fill[color=black](4+0, 0)circle(.5mm);
\fill[color=black](4+0.5, 0.75)circle(.5mm);
\fill[color=black](4+-0.5, 1)circle(.5mm);
\fill[color=black](6+0, 0)circle(.5mm);
\fill[color=black](6+0.5, 0.75)circle(.5mm);
\fill[color=black](6+-0.5, 1)circle(.5mm);
\draw[-](0, 0)--(2, 0);
\draw[-](2, 0)--(4, 0);
\draw[-](4, 0)--(6, 0);
\draw[-](0, 0)--(.5, 0.75);
\draw[-](0, 0)--(-.5, 1);
\draw[-](.5, 0.75)--(-.5, 1);
\draw[-](2+0, 0)--(2+.5, 0.75);
\draw[-](2+0, 0)--(2+-.5, 1);
\draw[-](2+.5, 0.75)--(2+-.5, 1);
\draw[-](4+0, 0)--(4+.5, 0.75);
\draw[-](4+0, 0)--(4+-.5, 1);
\draw[-](4+.5, 0.75)--(4+-.5, 1);
\draw[-](6+0, 0)--(6+.5, 0.75);
\draw[-](6+0, 0)--(6+-.5, 1);
\draw[-](6+.5, 0.75)--(6+-.5, 1);
\draw[-](0+.5, 0.75)--(6+.5, 0.75);
\path(+3, -1) node {The graph of $\Z\times_\Ic \Z/3\Z$, with $|\Ic|=2$};
\end{tikzpicture}
\end{align*}
\begin{align*}
\begin{tikzpicture}[ scale=1]
\path(-.5, 0) node {$\cdots$};
\path(6.5, 0) node {$\cdots$};
\fill[color=black](0, 0)circle(.5mm);
\fill[color=black](0.5, 0.75)circle(.5mm);
\fill[color=black](-0.5, 1)circle(.5mm);
\fill[color=black](2+0, 0)circle(.5mm);
\fill[color=black](2+0.5, 0.75)circle(.5mm);
\fill[color=black](2+-0.5, 1)circle(.5mm);
\fill[color=black](4+0, 0)circle(.5mm);
\fill[color=black](4+0.5, 0.75)circle(.5mm);
\fill[color=black](4+-0.5, 1)circle(.5mm);
\fill[color=black](6+0, 0)circle(.5mm);
\fill[color=black](6+0.5, 0.75)circle(.5mm);
\fill[color=black](6+-0.5, 1)circle(.5mm);
\draw[-](0, 0)--(2, 0);
\draw[-](2, 0)--(4, 0);
\draw[-](4, 0)--(6, 0);
\draw[-](0+.5, 0.75)--(2+.5, 0.75);
\draw[-](2+.5, 0.75)--(4+.5, 0.75);
\draw[-](4+.5, 0.75)--(6+.5, 0.75);
\draw[-](0-.5, 1)--(2, 1);
\draw[-](2-.5, 1)--(4, 1);
\draw[-](4-.5, 1)--(6-.5, 1);
\draw[-](0, 0)--(.5, 0.75);
\draw[-](0, 0)--(-.5, 1);
\draw[-](.5, 0.75)--(-.5, 1);
\draw[-](2+0, 0)--(2+.5, 0.75);
\draw[-](2+0, 0)--(2+-.5, 1);
\draw[-](2+.5, 0.75)--(2+-.5, 1);
\draw[-](4+0, 0)--(4+.5, 0.75);
\draw[-](4+0, 0)--(4+-.5, 1);
\draw[-](4+.5, 0.75)--(4+-.5, 1);
\draw[-](6+0, 0)--(6+.5, 0.75);
\draw[-](6+0, 0)--(6+-.5, 1);
\draw[-](6+.5, 0.75)--(6+-.5, 1);
\path(+3, -1) node {\emph{The graph of $\Z\times_\Ic \Z/3\Z$, with $|\Ic|=3$}};
\end{tikzpicture}
\end{align*}

 Under the representation $\ell^2(\Vc, m)\simeq \ell^2(\Vc_1, m_1)\otimes
\ell^2(\Vc_2, m_2)$, 
\begin{align}\label{e:deg_cp}
\deg_{\Gc}(\cdot)= \deg_{\Gc_1}(\cdot)\otimes \frac{1_\Ic(\cdot)}{m_2(\cdot)}
+ \frac{1}{m_1(\cdot)}\otimes \deg_{\Gc_2}(\cdot)
\end{align}
and
\begin{align}\label{e:rule}
\Delta_{\Gc, \theta} = \Delta_{\Gc_1, \theta_1} \otimes \frac{1_\Ic(\cdot)}{m_2(\cdot)} +
\frac{1}{m_1(\cdot)} \otimes \Delta_{\Gc_2,\theta_2}. 
\end{align}
If $m$ is non-trivial, we stress that the Laplacian obtained with our
product is usually not unitarily equivalent to the Laplacian obtained
with the Cartesian product. However, there is a potential $V:\Vc\to
\R$ such that $\Delta_{\Gc_1\times \Gc_2}$ is unitarily equivalent to
$\Delta_{\Gc_1\times_{\Vc_2} \Gc_2}+ V(\cdot)$, in $\ell^2(\Vc, m)$. 

\begin{definition}\label{d:dc}
Set $\Gc_1:=(\Ec_1, \Vc_1, m_1)$,
$\Gc_2:=(\Ec_2, \Vc_2, m_2)$, and $\Ic\subset \Vc_2$. We say that 
$\Gc= \Gc_1 \times_\Ic \Gc_2$ is a \emph{discrete cusp}
if the following hypotheses are satisfied:
\begin{itemize}
\item[(H1)] $m_1(x)$ tend to $0$ as $|x|\to \infty$, 
\item[(H2)] $\Gc_2$ is finite,
\item[(H3)]$\Delta_{\Gc_1, \theta_1}$ is bounded (or equivalently
  $\sup_{x\in \Vc_1}\deg_{\Gc_1}(x)<\infty$).
\end{itemize}
\end{definition}

We now motivate the choice of the above hypotheses by discussing the radius of
injectivity. We start by defining a different metric on $\Vc$, this choice is motivated by the works of 
\cite{CTT2} and \cite{MIT} but it needs a small adaptation for our purpose.

\begin{definition}
Given $\Gc:=(\Ec, \Vc, m)$,  
the \emph{weighted length} of an edge $(x,y)\in \Ec$ defined by:
\[L_\Gc\big((x,y) \big):=\sqrt{\frac{{\min\big(m(x),
    m(y)\big)}}{{\Ec(x,y)}}}.\]  
Given $x,y\in\Vc$, we define the \emph{weighted distance} from $x$ to
$y$ with respect to this length by:
\[\rho_{L_\Gc}(x,y):= \inf_{\gamma} \sum_{i=0}^{|\gamma|-1} L_\Gc\big( \gamma(i), \gamma({i+1})\big),\]
where $\gamma$ is a path joining $x$ to $y$ and with the convention that $\rho_{L_\Gc}(x,x):=0$ for all $x\in \Vc$.
\end{definition}

\begin{remark}
Since $\Gc$ is assumed connected, $\rho_{L_\Gc}$ is a metric on $\Vc$.
In fact $\rho_{L_\Gc}$  belongs to the class of \emph{intrinsic metrics}. We refer to  \cite{Ke} for a general definition, historical references,  properties, and applications. However, since Propositions \ref{p:radcusp} and \ref{p:radcusp2} do not hold in general with an arbitrary intrinsic metric, we stick to our specific choice of metric.
\end{remark}

We turn to the definitions of the girth and of the weighted radius of injectivity. This is essentially a weighted version of the standard ones, e.g, \cite{EGL}.
\begin{definition}
Given $\Gc:=(\Ec, \Vc, m)$, the \emph{girth} at $x\in \Vc$ of $\Gc$
w.r.t. the weighted length\ $L_\Gc$ is  
\[{\rm girth}(x):=\inf\{L_\Gc(\gamma), \gamma \mbox{ simple cycle of
  unweighted length } \geq 3 \mbox{ and containing } x\},\] 
where simple cycle means a closed walk with no repetitions of vertices and edges allowed, other than the repetition of the starting and ending vertex.
 We use the convention
 that the girth is $+\infty$ if there is no such cycle. 
\[{\rm girth}(\Gc):= \inf_{x\in \Vc}{\rm girth}(x).\]
The \emph{radius of injectivity} (at $x$) of $\Gc$ with respect to $L_\Gc$ is half the girth (at $x$). We denote the 
radius of injectivity by ${\rm rad}(\Gc)$ (at $x$ by ${\rm rad}(x)$ respectively)
\end{definition}
Note that with this definition, the radius of injectivity of a tree is
$+\infty$. 
 
\begin{proposition}\label{p:radcusp}
Given $\Gc_1:=(\Ec_1, \Vc_1, m_1)$ and
$\Gc_2:=(\Ec_2, \Vc_2, m_2)$ and $\Ic\subset \Vc_2$
Assume that $\Gc:= \Gc_1\times_\Ic \Gc_2$ is a discrete cusp. We have:
\begin{enumerate}[1)]
\item ${\rm rad}(\Gc_1)>0$.
\item If ${\rm rad}(\Gc_2)<\infty$, then ${\rm rad}(\Gc)=0$.
\end{enumerate}\end{proposition} 
\proof (1) Assume that ${\rm rad}(\Gc_1)=0$. Then for all $\varepsilon>0$, there is $x\sim y$ in $\Vc_1$ such that $L_{\Gc_1}\big((x,y)\big)<\varepsilon$. 
In particular, we have $\deg_{\Gc_1}(x)>\varepsilon^{-2}$ or  $\deg_{\Gc_1}(y)>\varepsilon^{-2}$. This is in contradiction with (H3).

(2) Since ${\rm rad}(\Gc_2)<\infty$, for all $x\in \Vc_1$, there is a
pure cycle contained in $\{x\}\times \Vc_2$. Moreover, for all $x\in
\Vc_1$ and $a\sim b$ in $\Vc_2$, since $\Ec(x,x)=0$, we have:
\begin{align*}
L_{\Gc_1\times_\Ic \Gc_2}\big(((x,a), (x,b))\big)&= \sqrt{m_1(x)} L_{\Gc_2}\big((a,b)\big)
\end{align*}
By (H1) we obtain that ${\rm rad}(\Gc)=0$. \qed

In contrast with this result we see that under the same hypotheses, the Cartesian product is not small at infinity. More precisely, we have:
\begin{proposition}\label{p:radcusp2}
Set $\Gc_1:=(\Ec_1, \Vc_1, m_1)$ and
$\Gc_2:=(\Ec_2, \Vc_2, m_2)$. Assume that (H1), (H2), and (H3) are satisfied. Then 
${\rm rad}(\Gc_1\times \Gc_2)>0$.
\end{proposition}
\proof 
Assume that ${\rm rad}(\Gc_1\times \Gc_2)=0$. For all $\varepsilon>0$, there are $x_1\sim y_1$ in $\Vc_1$ and $x_2\sim y_2$ in $\Vc_2$ such that 
\begin{align*}
\varepsilon > L_{\Gc_1\times \Gc_2}\big(((x_1,x_2), (x_1,y_2))\big)&=  L_{\Gc_2}\big((x_2,y_2)\big) 
\\
\mbox{ or } \varepsilon > L_{\Gc_1\times \Gc_2}\big(((x_1,x_2), (y_1,x_2))\big)&=  L_{\Gc_1}\big((x_1,y_1)\big).
\end{align*}
The first line is in  contradiction with (H2) and the second line with (H3). \qed

\subsection{Absence of essential spectrum}\label{s:destroy}

We have a first result of absence of essential spectrum. We refer to
\cite{CTT3} for   related  results based on the non-triviality of
$\Hol_\theta$ in the context of non-complete graphs. See also \cite{BGKLM} for similar ideas. 

\begin{proposition}\label{t:first}
Set $\Gc_1:=(\Ec_1, \Vc_1, m_1)$,
$\Gc_2:=(\Ec_2, \Vc_2, m_2)$, and $\Gc:=\Gc_1\times_\Ic \Gc_2$, with $|\Ic|>0$. Assume that  (H1), (H2), and 
 $\Hol_{\theta_2}\neq 0$ hold true. Then
$\Delta_{\Gc,   \theta}$ has a compact resolvent, and 
  \begin{align*}
\Nc_\lambda \left(m_1^{-1}(\cdot) \otimes
  \Delta_{\Gc_2,\theta_2}\right) \geq \Nc_\lambda(\Delta_{\Gc,\theta}), \mbox{ for all } \lambda\geq 0.
\end{align*}   
\end{proposition}

\proof 
 Note that 
\begin{align*}
\Delta_{\Gc,
  \theta}\geq \frac{1}{m_1(\cdot)} \otimes \Delta_{\Gc_2,\theta_2}
\end{align*}
in the form sense on $\Cc_c(\Vc)$. Since (H1) and (H2) hold, Lemma \ref{l:keylemm} ensures that $0$ is not in the spectrum of
$(\Delta_{\Gc_2, \theta_2})$. Hence  the spectrum of the
r.h.s.\ is purely discrete. By the min-max Principle, e.g., \cite{Go, RS}, 
the  operator $\Delta_{\Gc, \theta}$ has a compact resolvent.\qed 

\subsection{The asymptotic of the eigenvalues}\label{s:asymp}

From now on, we focus on the case when the graph is a discrete cusp and aim at a more precise result. 
To start off, we give the key-stone of our approach:

\begin{proposition}\label{p:keystone}
Set $\Gc_1:=(\Ec_1, \Vc_1, m_1)$,
$\Gc_2:=(\Ec_2, \Vc_2, m_2)$, and $\Ic\subset \Vc_2$ non-empty. Assume that  $\Gc:=\Gc_1\times_\Ic \Gc_2$ is a discrete cusp. 
We set
\begin{align} \label{e:M}
 M:=\sup_{x\in \Vc_1}\deg_{\Gc_1}(x)\times
\max_{y\in \Vc_2}(1/m_2(y))<\infty.
\end{align}
 We have:
\begin{align}\label{e:degimportant}
 \frac{1}{m_1(\cdot)}\otimes \deg_{\Gc_2}(\cdot)\leq \deg_{\Gc}(\cdot) \leq 
 \frac{1}{m_1(\cdot)}\otimes \deg_{\Gc_2}(\cdot) + M,
 \end{align}
\begin{align}\label{e:mainineq_cp}
\frac{1}{m_1(\cdot)} \otimes \Delta_{\Gc_2,\theta_2} \leq \Delta_{\Gc,
  \theta}\leq 2M+ \frac{1}{m_1(\cdot)} \otimes \Delta_{\Gc_2,\theta_2},
\end{align}
in the form sense on $\Cc_c(\Vc)$.
\end{proposition}
\proof Use \eqref{e:inedeg},  \eqref{e:deg_cp}, and \eqref{e:rule}. \qed

We work in the spirit of \cite{Go, BGK, BGKLM} and compare the Laplacian directly with the degree.

\begin{proposition}\label{p:firststep}
Set $\Gc_1:=(\Ec_1, \Vc_1, m_1)$,
$\Gc_2:=(\Ec_2, \Vc_2, m_2)$, and $\Ic\subset \Vc_2$ non-empty. Assume that  $\Gc:=\Gc_1\times_\Ic \Gc_2$ is a discrete cusp.  Set $M$ as in \eqref{e:M}. We have:
\begin{align}\label{e:formdeg}
\frac{\inf \sigma(\Delta_{\Gc_2,\theta_2})}{\max_{y\in \Vc_2}
 \deg_{\Gc_2}(y)} \left(\deg_{\Gc}(\cdot)-M\right)
\leq \Delta_{\Gc, 
  \theta}\leq 2M+2\deg_{\Gc}(\cdot),
\end{align}
in the form sense on $\Cc_c(\Vc)$.

Moreover, assuming that $\inf \sigma(\Delta_{\Gc_2,\theta_2})>0$, then
$\Dc(\Delta_{\Gc,
  \theta}^{1/2})=
\Dc\left(\deg_{\Gc}^{1/2}(\cdot)\right)$. Furthermore, since $\lim_{|x|\to
\infty} \deg_\Gc(x)=\infty$, $\Delta_{\Gc,\theta}$ has a compact resolvent and
\begin{align*}
0<\frac{\inf \sigma(\Delta_{\Gc_2,\theta_2})}{\max_{y\in \Vc_2}
 \deg_{\Gc_2}(y)} &\leq \liminf_{n\to \infty}
\frac{\lambda_n(\Delta_{\Gc,\theta})}{\lambda_n(\deg_\Gc(\cdot))}\leq
\limsup_{n\to \infty} 
\frac{\lambda_n(\Delta_{\Gc,\theta})}{\lambda_n(\deg_\Gc(\cdot))}
\leq 2.
\end{align*}
\end{proposition}
\proof Use \eqref{e:mainineq_cp} and  \eqref{e:inedeg} to get 
\begin{align*}
\frac{\inf \sigma(\Delta_{\Gc_2,\theta_2})}{\max_{y\in \Vc_2}
  \deg_{\Gc_2}(y)}\, \frac{1}{m_1(\cdot)} \otimes \deg_{\Gc_2}(\cdot)
\leq \Delta_{\Gc, 
  \theta}\leq 2M+ \frac{2}{m_1(\cdot)} \otimes \deg_{\Gc_2}(\cdot),
\end{align*}
Then apply  \eqref{e:degimportant} to obtain \eqref{e:formdeg}. Concerning the statement about the eigenvalue this follows from the standard
consequences of the min-max Principle, e.g., \cite{Go}. \qed

Here, trying to compare directly $\Delta_{\Gc, \theta}$ to $\deg_\Gc$ 
to get sharp results about eigenvalues is too optimistic because it is unclear 
how to obtain constants arbitrarily close to $1$ in front of
$\deg_\Gc$, as in \cite{Go, BGK}.  
To obtain some sharp asymptotics for the eigenvalues of $\Delta_{\Gc, \theta}$, as in \eqref{e:asymp_a}, 
we will use directly \eqref{e:mainineq_cp} and analyze very carefully the operator ${m_1^{-1}(\cdot)} \otimes \Delta_{\Gc_2,\theta_2}$. 

\begin{theorem}\label{t:main_cp}
Set $\Gc_1:=(\Ec_1, \Vc_1, m_1)$,
$\Gc_2:=(\Ec_2, \Vc_2, m_2)$, and $\Ic\subset \Vc_2$ non-empty. Assume that  $\Gc:=\Gc_1\times_\Ic \Gc_2$ is a discrete cusp. We  obtain that
 \begin{align}\label{e:main_cp}
\Dc(\Delta_{\Gc,
  \theta}^{1/2})= \Dc\left(m_1^{-1/2}(\cdot) \otimes
  \Delta_{\Gc_2,\theta_2}^{1/2} \right).
\end{align}
Moreover, we have:
\begin{enumerate}[1)]
\item $\Delta_{\Gc,
  \theta}$  has a compact resolvent if and only if
$\Hol_{\theta_2}\neq 0$.
\item If $\Hol_{\theta_2}\neq 0$, then
\[\Dc(\Delta_{\Gc,
  \theta}^{1/2})= \Dc\left(\deg_{\Gc}^{1/2}(\cdot)\right) \] 
and
\begin{align}\label{e:main_cp3}
\lim_{n\to \infty}\frac{\lambda_n \left(\Delta_{\Gc,
  \theta}\right)}{\lambda_n \left(m_1^{-1}(\cdot) \otimes
  \Delta_{\Gc_2,\theta_2}\right)}=1.
  \end{align}
Furthermore, setting $M$ as in \eqref{e:M}, 
\begin{align}\label{e:main_cp2}
\Nc_{\lambda-2M}\left(m_1^{-1}(\cdot)\otimes\Delta_{\Gc_2, \theta_2}\right)\leq
\Nc_\lambda(\Delta_{\Gc, \theta}) \leq
\Nc_\lambda\left(m_1^{-1}(\cdot)\otimes\Delta_{\Gc_2, \theta_2}\right),
\end{align}
for all $\lambda\geq 0$.
\end{enumerate}  
\end{theorem} 

\proof First note that \eqref{e:main_cp} follows directly from
\eqref{e:mainineq_cp}.  Denoting by $\{g_i\}_{i=1, .., 
|\Vc_2|}$ the  eigenfunctions associated to the eigenvalues
$\{\lambda_i\}_{i=1, .., |\Vc_2|}$ of $\Delta_{\Gc_2,\theta_2}$,  
where $\lambda_j\leq \lambda_{j+1}$, 
 we see that the eigenfunctions of
$m_1^{-1}(\cdot)\otimes \Delta_{\Gc_2}$ are given by
$\{\delta_{x}\otimes g_i\}$, where $x\in \Vc_1$ and $i=1, .., 
|\Vc_2|$. Then, using (H1), we observe that
\begin{align*}
\sigma\left(m_1^{-1}(\cdot)\otimes \Delta_{\Gc_2}\right) &= 
\overline{m_1^{-1}(\Vc_1)\times \{\lambda_1, \ldots, \lambda_{|\Vc_2|}\}}
= 
m_1^{-1}(\Vc_1)\times \{\lambda_1, \ldots, \lambda_{|\Vc_2|}\} .
\end{align*}
Besides, $0\in \sigma\left(m_1^{-1}(\cdot)\otimes
  \Delta_{\Gc_2}\right)$ if and only if $0$ is an eigenvalue of
$m_1^{-1}(\cdot)\otimes \Delta_{\Gc_2}$ of infinite multiplicity if
and only if  $\lambda_1=0$ if and only if $\Hol_{\theta_2}=
0$, by Lemma \ref{l:keylemm}. Moreover, recalling (H1), we see that all the eigenvalues of
$m_1^{-1}(\cdot)\otimes \Delta_{\Gc_2}$ which are not $0$ are of
finite multiplicity. Therefore, $m_1^{-1}(\cdot)\otimes
\Delta_{\Gc_2}$ has a compact resolvent if and only if $\Hol_{\theta_2} \neq 0$. 
Combining the latter and \eqref{e:mainineq_cp}, the min-max  Principle yields
the first point.

We turn to the second point and assume that $\Hol_{\theta_2}\neq
0$. The equality of the form-domains is given by \eqref{e:formdeg}.
Taking in account   \eqref{e:mainineq_cp}, the min-max Principle ensures
the asymptotic behavior of  $\lambda_n$ and the inequalities \eqref{e:main_cp2}. \qed

\begin{remark}\label{r:notdeg}
In the case when $\Hol_{\theta_2}=0$, for instance when $\theta_2=0$,
we see that the form-domain is $m_1^{-1/2}\otimes
P^{\perp}_{\ker(\Delta_{\Gc_2, \theta_2})}$. In particular, 
the form-domain is {\bf not} that of $\deg_\Gc(\cdot)$. Indeed
if the two form-domains are the same, the closed graph theorem yields
the existence of   
$c_1>0$ and $c_2>0$ so that
\[c_1 \deg_\Gc(\cdot) - c_2 \leq  m_1^{-1/2}\otimes
P^{\perp}_{\ker(\Delta_{\Gc_2, \theta_2})}, \]
in the form sense on $\Cc_c(\Vc)$. However, note that $0\in
\sigma_{\rm   ess}\left(m_1^{-1/2}\otimes P^{\perp}_{\ker(\Delta_{\Gc_2,
      \theta_2})}\right)$, whereas $\deg(\cdot)$ has a compact
resolvent. This is a contradiction with the min-max 
Principle. We obtain:
\begin{align*}
\Dc\left(\Delta_{\Gc,
  \theta}^{1/2}\right)= \Dc\left(\deg^{1/2}(\cdot)\right) &\Leftrightarrow \Hol_{\theta_2}\neq 0
  \\
  &\Leftrightarrow \Delta_{\Gc,
  \theta} \mbox{ has a compact resolvent.}
  \end{align*}
\end{remark}

In \eqref{e:main_cp3}, we exhibit the behaviour of the eigenvalues in terms of an explicit and computable mean. We now aim at comparing the asymptotic with that of the degree, as in \cite{Go, BGK}. The new phenomenon is that we are able to obtain a constant different from $1$ in the asymptotic.

 \begin{corollary}\label{c:referee}
 Let $\Gc_1:=(\Ec_1, \Vc_1, m_1)$,
$\Gc_2:=(\Ec_2, \Vc_2, m_2)$, and $\Ic\subset \Vc_2$ non-empty such that   $\Gc:=\Gc_1\times_\Ic \Gc_2$ is a discrete cusp.  Suppose that $\deg_{\Gc_2}$ is constant on $\Vc_2$ and take
  $\theta_2$ such that 
  $\Hol_{\theta_2}\neq 0$. Then, for all $a\in [1, +\infty[$, there exists $\widetilde \Gc_1:=(\widetilde\Ec_1,  \Vc_1, \tilde m_1)$ such that 
\begin{enumerate}[1)]
\item $\widetilde\Gc:=\widetilde\Gc_1\times_\Ic \Gc_2$ is a discrete cusp.
\item $\Ec_1$ and $\widetilde\Ec_1$ have the same zero set.
\item $\deg_{\widetilde \Gc_1}(x)\leq \deg_{\Gc_1}(x)$ for all $x\in \Vc_1$.
\item $\Delta_{\widetilde\Gc, \theta}$ is with compact resolvent, and
\begin{align}\label{e:asymp_a}
\lim_{\lambda\to \infty} \frac{\Nc_\lambda\left(\Delta_{\widetilde\Gc, \theta}\right)}{\Nc_\lambda\left(\deg_{\widetilde  \Gc}(\cdot)\right)}=a.
\end{align} 
\end{enumerate}
\end{corollary}
\proof  We choose $\widetilde m_1$ and $\widetilde \Ec_1$ later.
We denote by $\{\lambda_i\}_{i=1, \ldots, |\Vc_2|}$ the eigenvalues of $\Delta_{\Gc_2, \theta_2}$. Since $\Hol_{\theta_2}\neq 0$, we have $\lambda_i\neq 0$ for all $i=1, \ldots, |\Vc_2|$. This yields:
\begin{align*}
\Nc_\lambda\left(\frac{1}{\widetilde m_1}(\cdot)\otimes \Delta_{\Gc_2, \theta_2}\right) &= \left|\left\{(x, i), \frac{\lambda_i}{\widetilde m_1(x)} \leq \lambda\right\}\right|
=\sum_{i=1}^{|\Vc_2|} \left| \left(\frac{1}{\widetilde m_1}\right)^{[-1]} \left( \left[0, \frac{\lambda}{\lambda_i}\right]\right)\right|,
\end{align*}
where $[-1]$ denotes the reciprocal image.
On the other hand, 
\[\Nc_\lambda\left(\frac{1}{\widetilde m_1(\cdot)}\otimes \deg_{\Gc_2}\right)=
|\Vc_2|\times \left| \left(\frac{1}{\widetilde m_1}\right)^{[-1]} \left( \left[0,
      \frac{\lambda}{\deg_{\Gc_2}}\right]\right)\right|.\]
Moreover, from \eqref{e:degimportant} we get
\begin{align}\label{e:Ndeg}
\Nc_{\lambda-M}(\widetilde m_1^{-1}(\cdot)\otimes \deg_{\Gc_2})\leq 
\Nc_\lambda(\deg_{\widetilde\Gc}(\cdot)) \leq \Nc_\lambda(\widetilde m_1^{-1}(\cdot)\otimes \deg_{\Gc_2}),
\end{align}
for all $\lambda \geq 0$, where $M$ is given by \eqref{e:M}.
      
\noindent{\bf Step 1:} We first aim at $a=1$ in \eqref{e:asymp_a}. Thanks to Lemma \ref{l:asymp}, we choose $\widetilde m_1$ and $\widetilde \Ec_1$ such that the three first points are satisfied and
\[\left|\left\{x\in \Vc_1, \frac{1}{\widetilde m_1(x)} \leq \lambda\right\}\right|\sim \ln(\lambda), \quad \mbox{ as } \lambda\to \infty,\]
where $\sim$ stands for asymptotically equivalent. 
We obtain:
\begin{align}\label{e:lim0}
\frac{\Nc_\lambda\left(\frac{1}{\widetilde m_1(\cdot)}\otimes \Delta_{\Gc_2,
    \theta_2}\right)}{\Nc_\lambda\left(\frac{1}{\widetilde m_1(\cdot)}\otimes
    \deg_{\Gc_2}\right)}\sim \frac{\sum_{i=1}^{|\Vc_2|}(\ln(\lambda)-
  \ln(\lambda_i))}{|\Vc_2|(\ln(\lambda)- \ln(\deg_{\Gc_2}))}\to 1, \quad
\mbox{ as } \lambda\to \infty.
\end{align}
and for all $c\in\R$,
\begin{align}\nonumber
\Nc_{\lambda-c}\left(\frac{1}{\widetilde m_1(\cdot)}\otimes
    \deg_{\Gc_2}\right)& \sim |\Vc_2|\ln(\lambda-c) \sim |\Vc_2| \ln(\lambda)
    \\ \label{e:lim1}
    &\sim \Nc_\lambda\left(\frac{1}{\widetilde m_1(\cdot)}\otimes
    \deg_{\Gc_2}\right), \quad \mbox{ as } \lambda\to \infty.
\end{align} 
Combining the latter with \eqref{e:Ndeg}, we infer that for all $c\in \R$
\begin{align}
\Nc_{\lambda-c}\left(\frac{1}{\widetilde m_1(\cdot)}\otimes
    \deg_{\Gc_2}\right) \sim  \Nc_\lambda(\deg_{\widetilde\Gc}(\cdot)), \quad
\mbox{ as } \lambda\to \infty.
\end{align}
Using now \eqref{e:lim0}, this yields that for all $c\in \R$
\begin{align}
\Nc_{\lambda-c}\left(\frac{1}{\widetilde m_1(\cdot)}\otimes \Delta_{\Gc_2,
    \theta_2}\right) \sim  \Nc_\lambda(\deg_{\widetilde\Gc}(\cdot)), \quad
\mbox{ as } \lambda\to \infty.
\end{align}
Finally recalling \eqref{e:main_cp2}, we infer that
\begin{align*}
\Nc_\lambda\left(\Delta_{\widetilde \Gc, \theta}\right) \sim \Nc_\lambda\left(\deg_{\widetilde \Gc}(\cdot)\right), \quad
\mbox{ as } \lambda\to \infty.
\end{align*}
In other words, there are $\widetilde m_1$ and $\widetilde \Ec_1$ such that the three first points are satisfied and such that \eqref{e:asymp_a} is satisfied with $a=1$.

\noindent{\bf Step 2:} We turn to the case $a>1$ in \eqref{e:asymp_a}. Given $\alpha>0$,. Thanks to Lemma \ref{l:asymp}, we choose $\widetilde m_1$ and $\widetilde \Ec_1$ such that the three first points are satisfied and
\[\left|\left\{x\in \Vc_1, \frac{1}{\widetilde m_1(x)} \leq \lambda\right\}\right|\sim \lambda^{\alpha}, \quad \mbox{ as } \lambda\to \infty,\]
We obtain: 
\[\frac{\Nc_\lambda\left(\frac{1}{\widetilde m_1(\cdot)}\otimes \Delta_{\Gc_2,
    \theta_2}\right)}{\Nc_\lambda\left(\frac{1}{\widetilde m_1(\cdot)}\otimes
    \deg_{\Gc_2}\right)}\build{\sim}_{\lambda \to \infty}^{}
\frac{1}{|\Vc_2|}\sum_{i=1}^{|\Vc_2|} \left(\frac{ \deg_{\Gc_2}}{\lambda_i}\right)^\alpha=: F(\alpha).\]
First note that 
\[\lim_{\alpha\to 1^+} F(\alpha) = 1.\]
Next,  the sum of the eigenvalues (counted with multiplicity) of  $\Delta_{\Gc_2, \theta_2}$ is equal to 
$|\Vc_2|\deg_{\Gc_2}$. Therefore, there exists at least one eigenvalue $\lambda_i$, with $1 \leq i \leq |\Vc_2|$ so that $\deg_{\Gc_2}> \lambda_i$. In particular
\[\lim_{\alpha\to +\infty} F(\alpha) = +\infty.\]
Finally,  by continuity of $F$, we obtain that for all $a>1$ there is $\alpha>1$ such that $F(\alpha)= a$. To conclude, repeating the end of step 1, we obtain that for all  $a>1$, there   are $\widetilde m_1$ and $\widetilde \Ec_1$ such that the three first points are satisfied and such that \eqref{e:asymp_a} is satisfied.\qed

\begin{remark}
In \cite{Go, BGK}, the asymptotic in $\Nc_\lambda$ was not discussed
since the estimates that they obtain seem too weak to conclude. Being
able to compute $\Nc_\lambda$ in an explicit way, as in
\eqref{e:asymp_a}, is a new phenomenon. 
\end{remark}
We have used the following lemma:
\begin{lemma}\label{l:asymp}
Let $\Gc_1:=(\Ec_1, \Vc_1, m_1)$ be a graph satisfying (H1) and (H3) in Definition \ref{d:dc}   and let $f:[1, +\infty) \to [1, +\infty)$ be a continuous and strictly increasing function that tends to $+\infty$ at $+\infty$. There exists $\widetilde \Gc_1:=
(\widetilde \Ec_1, \Vc_1, \widetilde m_1)$ such that 
\begin{enumerate}[1)]
\item $\Ec$ and $\tilde\Ec$ have the same zero set.
\item (H1) and (H3) are satisfied for $\widetilde\Gc_1$.
\item $\deg_{\widetilde\Gc_1}(x)\leq \deg_{\Gc_1}(x)$ for all $x\in \Vc_1$.
\item We have:
\[\left|\left\{x\in \Vc_1, \frac{1}{\widetilde m_1(x)} \leq \lambda\right\}\right|\sim f(\lambda), \quad \mbox{ as } \lambda\to \infty,\]
 where $\sim$ stands for asymptotically equivalent.
\end{enumerate}
\end{lemma}
\proof 
Without any loss of generality, one may suppose that $f(1)=1$. Let $\phi:\N^* \to \Vc_1$ be a bijection. Set:
\[\tilde m_1(\phi(n)):= \frac{1}{f^{[-1]}(n)},\]
where $[-1]$ denotes the reciprocal image. Note that (H1) is satisfied. Moreover, 
\begin{align*}
\left|\left\{x\in \Vc_1, \frac{1}{\widetilde m_1(x)} \leq \lambda\right\}\right|= \left|\left\{n\in \N^*, n \leq f(\lambda)\right\}\right| = \lfloor f(\lambda)\rfloor +1 \sim f(\lambda),
\end{align*}
as $\lambda\to\infty$. Finally, we set:
\[\widetilde\Ec_1(x,y):= \Ec_1(x,y)\frac{\min(\tilde m_1(x), \tilde m_1(y))}{\max(m_1(x), m_1(y))}.\]
The first point is clear. For (H3), note that $\deg_{\widetilde\Gc_1}(x)\leq \deg_{\Gc_1}(x)$ for all $x\in \Vc_1$. \qed

We end this section by proving the results stated in the introduction.
\proof[Proof of Theorem \ref{t:intro}]
Let us consider $\Gc_1:=(\Ec_1, \Vc_1, m_1)$, where 
\[\Vc_1:=\N, \quad m_1(n):=\exp(-n), \mbox{ and } \Ec_1(n, n+1):= \exp(-(2n+1)/2),\]
for all $n\in \N$ and $\Gc_2:=(\Ec_2, \Vc_2, 1)$ a simple connected finite graph such that $|\Vc_2|=n$. 
Set $\Gc:=\Gc_1\times_{\Vc_2} \Gc_2$, $\theta_1:=0$ and $\theta_2$
such that $\Hol_{\theta_2}\neq 0$. 

In the spirit of \cite{GoMo}, we denote by  $P^{{\rm
    le}}_\kappa$  the projection on $\ker(\Delta_{\Gc_2, \kappa 
  \theta_2})$ and by $P^{{\rm
    he}}_\kappa$ is the projection on $\ker(\Delta_{\Gc_2,   \kappa
  \theta_2})^\perp$.  Here \emph{le} stands for \emph{low energy} and
\emph{he} for  \emph{high energy}.

We have that $\Delta_{\Gc, \kappa \theta}:= \Delta_{\Gc, \kappa \theta}^{{\rm
    le}}\oplus \Delta_{\Gc, \kappa \theta}^{{\rm he}}$, where
\begin{align*}
\Delta_{\Gc, \kappa \theta}^{{\rm le}} := \Delta_{\Gc_1, 0} \otimes P^{{\rm le}}_\kappa,  
\end{align*}
on $(1 \otimes P^{{\rm le}}_\kappa) \ell^2(\Vc, m)$, and 
\begin{align*}
\Delta_{\Gc, \kappa \theta}^{{\rm he}} := \Delta_{\Gc_1, 0} \otimes P^{{\rm he}}_\kappa+ 
\frac{1}{m_1(\cdot)} \otimes P^{{\rm he}}_\kappa \Delta_{\Gc_2, \kappa \theta_2}, 
\end{align*}
on $(1 \otimes P^{{\rm he}}_\kappa) \ell^2(\Vc, m)$.

By Lemma \ref{l:keylemm}, Corollary \ref{c:coupling}, and Remark \ref{r:notdeg}, there exists 
$\nu>0$ such that 
\begin{align*}
P^{{\rm le}}_\kappa=0 \quad &\Leftrightarrow \quad\Hol_{{\kappa\theta_2}}\neq 0
\\
& \Leftrightarrow \quad \kappa \neq 0 \mbox{ in }{\R/\nu\Z} \quad
\Leftrightarrow \quad \Dc\left(\Delta_{\Gc,\kappa \theta}^{1/2}\right) =
\Dc\left(\deg_\Gc^{1/2}(\cdot)\right).
\end{align*}
The proof of Theorem \ref{t:main_cp} gives the first point. Assume
that $\kappa \in {\R/\nu\Z}$. Let $U:\ell^2(\N, m_1)\to \ell^2(\N, 1)$
be the unitary map given by $Uf(n):= \sqrt{m_1(n)} f(n)$. We see that:
\begin{align*}
U\Delta_{\Gc, \kappa
  \theta}^{{\rm le}}U^{-1}=\Delta_{\N, 0}+ (e^{-1/2}-1)\delta_0 +
e^{1/2}+e^{-1/2}-2 \mbox{ in } \ell^2(\N),
\end{align*}
where $\Delta_{\N, 0}$ is related to the simple graph of
  $\N$. By using for instance some Jacobi matrices techniques, it is well-known that the  essential spectrum of $\Delta_{\Gc, \kappa
    \theta}^{{\rm le}}$ is purely  
 absolutely continuous and equal to 
 \[\sigma_{\rm ac}(\Delta_{\Gc, \kappa
    \theta}^{{\rm le}})= [e^{1/2}+e^{-1/2}-2,
 e^{1/2}+e^{-1/2}+2],\]
 with multiplicity one, e.g., \cite{Wei}. It has a unique
 eigenvalue and it is negative.  
 
 We turn to the high energy part. 
Denote by $\{\lambda_i\}_{i=1, \ldots, n}$, with $\lambda_i\leq
\lambda_{i+1}$, the eigenvalues of 
$\Delta_{\Gc_2, \kappa\theta_2}$. Recall that $\lambda_1 =0 $ due to
the fact that $\Hol_{{\kappa\theta_2}}= 0$. 
By  \eqref{e:mainineq_cp},
\begin{align*}
\frac{1}{m_1(\cdot)} \otimes \Delta_{\Gc_2,\kappa\theta_2}P^{{\rm
    he}}_\kappa \leq \Delta_{\Gc, 
 \kappa  \theta} (1\otimes P^{{\rm he}}_\kappa)\leq 2M+ \frac{1}{m_1(\cdot)}
\otimes \Delta_{\Gc_2,\kappa\theta_2}P^{{\rm he}}_\kappa.
\end{align*}
Hence, $\Delta_{\Gc, 
  \kappa\theta} (1\otimes P^{{\rm he}}_\kappa)$ has a compact resolvent and
\begin{align*}
\Nc_{\lambda-2M}\left(m_1^{-1}(\cdot)\otimes\Delta_{\Gc_2,
    \kappa\theta_2}P^{{\rm he}}_\kappa\right)&\leq 
\Nc_\lambda\left(\Delta_{\Gc, \kappa\theta} (1\otimes P^{{\rm
      he}}_\kappa)\right)  
\\
&\quad \quad\quad\quad\quad\leq
\Nc_\lambda\left(m_1^{-1}(\cdot)\otimes\Delta_{\Gc_2, \kappa\theta_2}P^{{\rm
      he}}_\kappa\right), 
\end{align*}
for all $\lambda\geq 0$.  Finally: 
\[\frac{\Nc_\lambda(\frac{1}{m_1(\cdot)}\otimes \Delta_{\Gc_2,
   \kappa  \theta_2}P^{{\rm he}}_\kappa)}{\Nc_\lambda\left(\frac{1}{m_1(\cdot)}\otimes
    \deg_{\Gc_2}\right)}\sim \frac{\sum_{i=2}^{n}\ln(\lambda)-
  \ln(\lambda_i)}{n(\ln(\lambda)- \ln(\deg_{\Gc_2}))}\to \frac{n-1}{n},
\mbox{ as } \lambda\to \infty.\]
We conclude with \eqref{e:lim1} for $a=1$. \qed

\end{document}